\documentclass[preprint]{elsarticle}

\makeatletter
\def\ps@pprintTitle{%
  \let\@oddhead\@empty
  \let\@evenhead\@empty
  \def\@oddfoot{\reset@font\hfil\thepage\hfil}
  \let\@evenfoot\@oddfoot
}
\makeatother

\usepackage{amssymb}
\usepackage{siunitx}
\usepackage{amsmath}
\usepackage{seqsplit}
\usepackage{graphicx}
\usepackage{amsmath}
\usepackage{amsthm}
\usepackage{amssymb}
\usepackage{multicol}
\usepackage{wrapfig}
\usepackage{fontawesome5}
\usepackage{multirow}
\usepackage{color}
\usepackage{hyperref}
\usepackage{booktabs}

\usepackage{wasysym}
\usepackage{url}  
\usepackage[figuresright]{rotating}
\usepackage{mathtools}
\usepackage{booktabs}
\usepackage{bigdelim}
\usepackage{enumitem}
\usepackage{algorithm}
\usepackage{csquotes}
\usepackage[noend]{algpseudocode}

\newcommand{\black}[1]{{\color{black}#1}}

\newtheorem{theorem}{Theorem}[section] 
\newtheorem{lemma}[theorem]{Lemma}     
\newtheorem{proposition}[theorem]{Proposition}
\newtheorem{corollary}[theorem]{Corollary}
\theoremstyle{definition} 
\newtheorem{definition}{Definition}[section]

\theoremstyle{remark} 
\journal{Operations Research Letters}

\begin{document}

\begin{frontmatter}



\title{\black{Parallel} Graver Basis Extraction for Nonlinear Integer Optimization}

%

\author[ucas]{Wenbo Liu} 
\author[sribd,cuhksz]{Akang Wang\corref{cor1}}
\author[ucas]{Wenguo Yang}

\cortext[cor1]{Corresponding author: Akang Wang \textless email: wangakang@sribd.cn\textgreater.}

\affiliation[ucas]{organization={University of Chinese Academy of Sciences},
        country={China}}


 \affiliation[sribd]{organization={Shenzhen International Center for Industrial and Applied Mathematics, Shenzhen Research Institute of Big Data},
             country={China}}

\affiliation[cuhksz]{organization={The Chinese University of Hong Kong},
            city={Shenzhen},
            country={China}}
            
\begin{abstract}

\black{The augmentation scheme provides a nontraditional approach to nonlinear integer programming by iteratively refining incumbent solutions along objective-improving directions from the Graver basis. 
Its main computational bottleneck, however, lies in the practical difficulty of accessing such directions. 
To address this challenge, we develop a massively parallel heuristic for approximating Graver basis, extracting promising directions by optimizing nonconvex continuous problems using parallelizable first-order methods. 
Experiments on QPLIB and MINLPLib instances show that our method achieves comparable performance to advanced solvers.}
\end{abstract}


   
    


\end{frontmatter}



\section{Introduction}
\label{sec:introduction}
\black{We consider a nonlinear integer program of the form}
\begin{subequations}
\label{eq:INLP}
\begin{alignat}{2}
    \underset{x\in\mathbb Z^n}{\min} & \quad f(x)\\
    \;\text{s.t.} & \quad Ax=b\label{eq:constr_linear}\\
                  & \quad l\leq x\leq u, \label{eq:constr_bounds}
\end{alignat}
\end{subequations}
where $f:\mathbb R^n\rightarrow\mathbb R$ is a real-valued function, and $A\in\mathbb Z^{m  \times n}, b\in\mathbb Z^m, l, u\in\mathbb{Z}^n$ define linear constraints and variable bounds.
Let $S\coloneqq\{x\in\mathbb Z^n:Ax=b,l\leq x\leq u\}$ denote the feasible region\black{, and let $\mathrm{Ker}_\mathbb Z(A)\coloneqq \{g\in\mathbb Z^n:Ag=0\}$ denote the \textit{integer kernel} of $A$}.
Such problems arise in real-world applications like portfolio selection~\cite{markowitz1967portfolio} and resource allocation~\cite{li2006nonlinear}, but are computationally challenging due to their \( \mathcal{NP} \)-hard complexity. General-purpose solvers rely on branch-and-cut methods, which are expensive due to systematic solution enumeration.

\black{A nontraditional approach for tackling nonlinear integer programs of form~(\ref{eq:INLP}) involves augmenting via the Graver basis. Given $x,y\in\mathbb R^n$, we say that \( x \sqsubseteq y \) if \( x_i y_i \geq 0 \) and \( |x_i| \leq |y_i| \) for all \( i = 1, 2, \dots, n \), and the Graver basis $\mathcal G(A)$ consists of all \( \sqsubseteq \)-minimal elements in  \( \mathrm{Ker}_{\mathbb{Z}}(A) \setminus \{0\} \).
The augmentation scheme begins with an initial solution $\bar x$ and iteratively refines it by moving along feasible, objective-reducing directions drawn from $\mathcal G(A)$.}
The Graver augmentation procedure is shown to achieve an optimal solution within polynomially many steps for separable convex objectives~\cite{hemmecke2011polynomial}.
Furthermore, the Graver basis is applicable to various nonlinear integer programs including convex integer maximization~\cite{de2009convex} and minimization of quadratic or higher degree polynomial functions lying in suitable cones~\cite{lee2012quadratic}.

\black{
An immediate question is how to access the Graver basis for practical augmentation. 
This difficulty arises from its intrinsic intractability: the size of $\mathcal G(A)$ may grow exponentially, and deciding membership in $\mathcal G(A)$ is $\mathcal{NP}$-complete~\cite{Henk1996OnHB}. 
Exact methods, including Pottier’s algorithm~\cite{pottier1996euclidean} and the project-and-lift algorithm~\cite{hemmecke2002computation}, are inherently sequential and thus poorly scalable. 
Recent works~\cite{alghassi2019graver,karahalios2024quantum} employs quantum annealing to generate diverse elements of $\mathrm{Ker}_{\mathbb Z}(A)$, but still depends on Pottier’s procedure to ensure $\sqsubseteq$-minimality.}

\black{In this work, we shift from exact computation to the approximation of Graver basis.}
This is performed by optimizing unconstrained, non-convex continuous problems. 
The objective function is carefully designed through an isomorphism to $\mathrm{Ker}_\mathbb Z(A)$, and first-order methods are employed for parallel optimization, initiated from diversified starting points. 
During this process, near-minimal elements of $\mathrm{Ker}_\mathbb{Z}(A)$ are progressively extracted to form this approximate basis.
This set of directions is then applied in multiple augmentations from diverse initial solutions, each striving for optimality. 
We refer to this entire process as \underline{M}ulti-start \underline{A}ugmentation via \underline{P}arallel \underline{E}xtraction (MAPE), a \black{massively parallel} algorithm for solving nonlinear integer programs.
MAPE exhibits three notable characteristics: (i)~it does not rely on general-purpose solvers, while ensuring feasible returned solutions;  
(ii)~it is highly parallelizable with GPU acceleration;
(iii)~it is reusable, as the parallel extraction needs to be performed only once to handle instances with identical constraint matrices but varying objectives and right-hand sides.

The distinct contribution of our work can be summarized as follows:
\begin{itemize}
    \item We \black{extract diverse Graver basis elements} by optimizing unconstrained problems, demonstrating how continuous techniques can be applied to discrete problems.

    \item We introduce MAPE, a \black{massively parallel} algorithm for solving nonlinear integer programs, \black{building upon the partial Graver basis.}

    \item We conduct extensive computational experiments on \black {MINLPLib and} QPLIB instances, with results showing that our approach performs comparably to state-of-the-art general-purpose solvers.
\end{itemize}

\section{The Graver Augmentation}
\label{sec:preliminary}
\begin{definition}
\label{def:test set}
A set $\mathcal T \subseteq \mathbb Z^n$ is called a \textit{test set} (or optimality certificate) for model~(\ref{eq:INLP}) if, for any non-optimal $\bar x \in S$, there exists $g \in \mathcal T$ and $\lambda \in \mathbb Z_+$ such that: (i)~$\bar x + \lambda g \in S$ and (ii)~$f(\bar x + \lambda g) < f(\bar x)$.
\end{definition}

Given an instance of model~(\ref{eq:INLP}) and a set $\mathcal T$ of directions, the augmentation algorithm begins with an initial feasible solution and repeatedly selects the best augmenting direction from $\mathcal{T}$ to improve the incumbent, until no such a direction exists.
We outline this augmentation procedure in Algorithm~\ref{alg:augmentation}.

\begin{algorithm}
\caption{Augmentation Scheme}
\label{alg:augmentation}
\begin{algorithmic}[1] 
\Statex \textbf{INPUT:} An optimization problem in model~(\ref{eq:INLP}), a finite  set $\mathcal T$ and an initial feasible solution $\bar x\in S$
\Statex \textbf{OUTPUT:} {A solution} $x^*\in S$.
\While {$\mathcal A\coloneqq\{(g,\lambda):g\in\mathcal T,\lambda\in\mathbb Z_+,\bar x+\lambda g\in S,f(\bar x+\lambda g)<f(\bar x)\}\neq \emptyset$}
\State $(\bar g,\bar\lambda)\gets\underset{(g,\lambda)\in\mathcal A} {\arg\min}\;f(\bar x+\lambda g)$
\State $\bar x\gets\bar x+\bar\lambda \bar g$.
\EndWhile
\State \textbf{return} $\bar x$
\end{algorithmic}
\end{algorithm}

\black{It can be easily shown that Algorithm~\ref{alg:augmentation} always returns an optimal solution to model~(\ref{eq:INLP}) whenever $\mathcal T$ is a test set.}
 The following proposition characterizes an important special case where $\mathcal{G}(A)$ forms a test set.

\begin{proposition}[\cite{de2012algebraic}]
\label{prop:graver_is_test_set}
The Graver basis $\mathcal{G}(A)$ forms a test set for model~(\ref{eq:INLP}) when the objective function $f: \mathbb{R}^n \to \mathbb{R}$ is separable convex, i.e., when $f(x) = \sum_{j=1}^n f_j(x_j)$ where each $f_j: \mathbb{R} \to \mathbb{R}$ is convex.
\end{proposition}

Proposition~\ref{prop:graver_is_test_set} establishes that for separable convex objectives, the Graver augmentation is guaranteed to reach an optimal solution to model~(\ref{eq:INLP}). Moreover, as shown in~\cite{hemmecke2011polynomial}, the number of required augmenting steps is polynomially bounded in the input size.

\section{MAPE: Multi-start Augmentation via Parallel Extraction}
\label{sec:extraction}
\black{To facilitate practical augmentation, we develop a parallel approach for extracting Graver basis elements, anchored in four fundamental conditions:}

\begin{itemize}[left=2em]
    \item[(C1)]\refstepcounter{enumi}\label{condition:kernel}Kernel condition: $Ag = 0$;
    \item[(C2)]\refstepcounter{enumi}\label{condition:integral}Integrality: $g \in \mathbb{Z}^n$;
    \item[(C3)]\refstepcounter{enumi}\label{condition:minimal}Minimality: $g$ is $\sqsubseteq$-minimal;
    \item[(C4)]\refstepcounter{enumi}\label{condition:bounds}Bound consistency: $l - u \leq g \leq u - l$.
\end{itemize}

Conditions (C\ref{condition:kernel}) and (C\ref{condition:integral}) characterize elements of $\mathrm{Ker}_{\mathbb{Z}}(A)$, and (C\ref{condition:minimal}) ensures non-redundant directions that define $\mathcal{G}(A)$. Notably, $\mathcal G(A)$ depends solely on the linear constraints~(\ref{eq:constr_linear}), whereas the variable bounds~(\ref{eq:constr_bounds}) admit further pruning of invalid directions through (C\ref{condition:bounds}).
Existing approaches like~\cite{alghassi2019graver} strictly enforce (C\ref{condition:integral}) and~(C\ref{condition:bounds}) via quantum bit encoding, followed by quantum annealing to minimize $\|Ag\|_2^2$ for kernel condition, and finally apply the Pottier's algorithm for minimality.
In contrast, we propose an alternative scheme that strictly enforces conditions (C\ref{condition:kernel}) and (C\ref{condition:integral}), and addresses both (C\ref{condition:minimal}) and (C\ref{condition:bounds}) through a \black{unified criterion}. 
Our approach first characterizes the integer kernel $\mathrm{Ker}_{\mathbb{Z}}(A)$ through its lattice basis representation, which inherently satisfies both the kernel condition and integrality requirement. 
We then develop a \black{parallel extraction} method that minimizes an unconstrained, non-convex continuous problem designed to promote direction minimality and bound consistency, progressively collecting near-minimal elements during the optimization.
\black{The resulting set of directions serves as an approximate Graver basis, enabling practical multi-start augmentation for nonlinear integer programs,} which we designate as MAPE.

\subsection{Characterization of $\mathrm{Ker}_\mathbb Z(A)$}
\label{sec:xz_scheme}
\black{For any full-column-rank matrix $B=(B_1,...,B_d)\in\mathbb R^{n\times d}$, denote by $\mathcal L(B)\coloneqq\left\{\sum\limits_{i=1}^d z_iB_i:z\in\mathbb Z^d\right\}$ the lattice generated by $B$.
We now construct a basis for $\mathrm{Ker}_\mathbb Z(A)$ of dimension $d\coloneqq n-m$.}

\begin{lemma}[\cite{kannan1979polynomial}]
\label{lem:hnf}
For $A\in\mathbb Z^{m\times n}$ with full row rank, there exists a unimodular matrix $K \coloneqq (D|B) \in\mathbb Z^{n \times (m+d)}$ such that
$AK=(H|0)$ is of Hermite normal form.
\end{lemma}


\black{The matrix $K$ can be computed in polynomial time~\cite{kannan1979polynomial}.
Its last $d$ columns, denoted by $B \coloneqq (B_1,\ldots,B_d)$, form a basis of $\mathrm{Ker}_{\mathbb Z}(A)$.}

\begin{theorem}
\label{thm:iso}
$\mathrm{Ker}_\mathbb Z(A)=\mathcal L(B).$
\end{theorem}

\begin{proof}
For any $Bz \in \mathcal{L}(B)$ where $z \in \mathbb{Z}^d$, Lemma~\ref{lem:hnf} gives $AB = 0$, hence $ABz = 0$. Since $B \in \mathbb{Z}^{n \times d}$ by Lemma~\ref{lem:hnf}, we have $Bz \in \mathbb{Z}^n$. Thus, $Bz \in \mathrm{Ker}_{\mathbb{Z}}(A)$.

For any $g \in \mathrm{Ker}_{\mathbb{Z}}(A)$, consider $K^{-1}g \in \mathbb{Z}^n$ (since $K \in \mathbb{Z}^{n \times (m+d)}$ is unimodular and $g \in \mathbb{Z}^n$). Decompose:
\[
K^{-1}g = \binom{w}{z} \quad \text{with } w \in \mathbb{Z}^m, z \in \mathbb{Z}^d
\]
Then:
\[
0 = Ag = AKK^{-1}g = (H|0)\binom{w}{z} = Hw
\]
The invertibility of $H$ implies $w = 0$. Therefore:
\[
g = K\binom{0}{z} = Bz \in \mathcal{L}(B).
\]
Therefore, $\mathrm{Ker}_\mathbb Z(A)=\mathcal L(B).$
\end{proof}

Theorem~\ref{thm:iso} indicates an isomorphism between $\mathrm{Ker}_\mathbb Z(A)$ and $\mathbb Z^d$ via $g=Bz$, and the inverse is given by $z=(B^\top B)^{-1}B^\top g$.
Therefore, every element $g \in \mathrm{Ker}_\mathbb Z(A)$ can be obtained via its integer coordinate in the space of $\mathbb Z^d$.

\subsection{Parallel Extraction of the Graver Basis}
\label{sec:parallel extraction}


Our approach enforces condition~(C\ref{condition:kernel}) through the parameterization $g = Bz$, where $B$ denotes a basis for the integer kernel $\mathrm{Ker}_{\mathbb{Z}}(A)$  constructed following ~\cite{kannan1979polynomial}. 
When necessary, the basis may be further refined using the LLL algorithm (see~\ref{apdx:lll} for details).
Theorem~\ref{thm:iso} establishes that condition~(C\ref{condition:integral}) can be satisfied through simple rounding of $z$. 
However, condition~(C\ref{condition:minimal}) presents greater challenges. 
To address this, we develop a unified criterion that simultaneously incorporates conditions~(C\ref{condition:minimal}) and~(C\ref{condition:bounds}) by identifying \emph{short} directions. This leads to the optimization problem:
\begin{equation}
\label{eq:SVP}    
\underset{z\in\mathbb Z^{d}\backslash\{0\}}{\min} \Vert Bz\Vert.
\end{equation}

Problem~(\ref{eq:SVP}) represents the shortest vector problem, which is known to be $\mathcal{NP}$-hard~\cite{svp}. 
However, instead of exclusively seeking the shortest vector, our goal is to generate a diverse collection of short vectors through a highly parallelizable implementation.
To this end, we optimize a non-convex continuous surrogate of Problem~(\ref{eq:SVP}), reformulated as the objective function in~(\ref{eq:loss}). 
In this formulation, we employ the $\ell_1$-norm, motivated by the well-established success of LASSO~\cite{lasso} in promoting sparsity. 
The model includes two additional components: (i) a penalty term that drives $z$ toward integer values in $\mathbb{Z}^d$, and (ii) a term that penalizes solutions near the origin through an inverse $\ell_\infty$-norm, with $\lambda_1$ and $\lambda_2$ as tuning parameters.
\black{The \textit{surrogate problem} is formulated as}:
\begin{equation}
\label{eq:loss}
    \underset{z\in\mathbb R^d}{\min} \; \Phi(z)\coloneqq \Vert Bz\Vert_1 + \lambda_1\cdot\sum\limits_{i=1}^d (z_i-\lfloor z_i\rfloor)(\lceil z_i\rceil -z_i)+\lambda_2\cdot \max \left\{\frac{1}{\Vert z\Vert_\infty}-1,0 \right\}.
\end{equation}

Given Problem~(\ref{eq:loss}), our parallel extraction method performs simultaneous optimization from multiple diverse initial points through three main phases. 
First, we generate $N$ uniformly distributed random points $g \sim \text{Uniform}(l-u, u-l)$ to ensure diversity in initialization. These points are projected via the transformation $z \coloneqq (B^\top B)^{-1}B^\top g$ to obtain valid starting points for model~(\ref{eq:loss}).
The optimization then proceeds using first-order methods, particularly Adam~\cite{kingma2014adam}, selected for their computational efficiency and native GPU parallelization support. During optimization, intermediate solutions are periodically processed through rounding and mapping to the integer kernel via $B\lceil z \rfloor$, yielding a comprehensive collection of candidate directions in $\mathrm{Ker}_{\mathbb{Z}}(A) \setminus \{0\}$.
This parallel extraction framework is formally described in Algorithm~\ref{alg:extracting_gb}.

\begin{algorithm}
\caption{Parallel Extraction}
\label{alg:extracting_gb}
\begin{algorithmic}[1] 
\Statex \textbf{INPUT:} $A\in \mathbb Z^{m\times n}$, number of epochs $T$, number of initial points $N$
\Statex \textbf{OUTPUT:} the collected direction set $\mathcal G$
\State $\mathcal G\gets \emptyset$
\State Compute the HNF of $A$ with $AK=(H|0)$
\State $B\gets K_{:,m+1:n}$
\For {$i=1,2,...,N$} \hfill \# in parallel
    \State Generate a sample $g \sim \text{Uniform}(l-u, u-l)$
    \State $z^0 \gets (B^\top B)^{-1}B^\top g$
\For {$t = 1, 2, ..., T$}
    \State $z^{t} \gets update(z^{t-1}, \nabla \Phi(z^{t-1}))$\hfill \# a single-step update
    \If {$B \lceil z^t \rfloor\in[l-u,u-l]$}  
    \State $\mathcal G \gets \mathcal G\cup\{B \lceil z^t \rfloor\}$
    \EndIf
\EndFor
\EndFor
\State \textbf{return} $\mathcal G$
\end{algorithmic}
\end{algorithm}

\subsection{Multi-start Augmentation via Parallel Extraction}
\label{sec:MAPE}
Section~\ref{sec:parallel extraction} provides a \black{parallel} approach for Graver basis extraction.
Building upon this, we introduce the MAPE algorithm for solving nonlinear integer programs. In particular, as the parallel extraction only yields a \black{partial} Graver basis, MAPE compensates by performing multiple augmentations in parallel.
\black{Initial feasible solutions to model~(\ref{eq:INLP}) are obtained by solving the following problem}:
\begin{equation}
\label{eq:init_sol}
    \underset{x\in\mathbb [l,u]}{\min} \Vert Ax-b\Vert_2^2 + \lambda_3\cdot\sum\limits_{i=1}^n(x_i-\lfloor x_i\rfloor)(\lceil x_i\rceil - x_i)
\end{equation}
where $\lambda_3$ is a hyper-parameter that balances integrality against satisfaction of linear constraints, and any optimal solution to Problem~(\ref{eq:init_sol}) is feasible for model~(\ref{eq:INLP}).
\black{Starting from diversified initial points, Problem~(\ref{eq:init_sol}) can be optimized in parallel using Adam}, as discussed in Section~\ref{sec:parallel extraction}.
Let $S_0\subseteq S$ denote the set of feasible solutions to model~(\ref{eq:INLP}).
Beginning from each $\bar x\in S_0$, MAPE performs $|S_0|$ parallel augmentations via the extracted  partial Graver basis and finally returns the best solution found.

\black{

\subsection{Extension to Cases with Linear Inequalities}
\label{apdx:lin_ineq}

Consider the presence of additional linear inequality constraints 
$Gx \ge h$, where $G \in \mathbb{Z}^{m' \times n}$ and 
$h \in \mathbb{Z}^{m'}$, in Problem~(\ref{eq:INLP}). 
Introducing slack variables $s \in \mathbb{Z}_+^{m'}$, 
the problem can be reformulated as
\begin{equation}
\label{eq:NLIP_ineq_slack}
\begin{aligned}
&\underset{x\in\mathbb Z^n,s\in\mathbb Z^{m'}}{\min}&&f(x)\\
&\quad \; \text{s.t.}&&\begin{pmatrix}
    A&0\\G&-I_{m'}
\end{pmatrix}\begin{pmatrix}
    x\\s
\end{pmatrix}=\begin{pmatrix}
    b\\h
\end{pmatrix}\\
& &&l\leq x\leq u,\;s\geq 0
\end{aligned}
\end{equation}
To characterize the integer kernel of the constraint matrix, we have the following result.
\begin{corollary}
\begin{equation*}
   \mathrm{Ker}_\mathbb Z\left(\begin{pmatrix}
    A&0\\G&-I_{m'}
\end{pmatrix}\right)=\mathcal L\left(\binom{B}{GB}\right), 
\end{equation*}where $B\in\mathbb Z^{n\times d}$ is a basis of $\mathrm{Ker}_\mathbb Z(A)$.
\end{corollary}

Accordingly, MAPE searches for \emph{short} directions in 
$\mathrm{Ker}_{\mathbb Z}\!\left(
\begin{pmatrix}
A & 0 \\
G & -I_{m'}
\end{pmatrix}
\right)$.
The corresponding surrogate problem becomes
\begin{equation}
\label{eq:surrogate_ineq}
\begin{aligned}
    &\underset{z\in\mathbb R^d}{\min}\;\left\Vert \binom{B}{GB}z\right\Vert_1+\lambda_1\sum\limits_{i=1}^d(z_i-\lfloor z_i\rfloor)(\lceil z_i\rceil-z_i)+\lambda_2\cdot \max\left\{\frac1{\Vert z\Vert_\infty}-1,0\right\}
\end{aligned}
\end{equation}
Compared with the standard surrogate problem~(\ref{eq:loss}), 
the objective in~(\ref{eq:surrogate_ineq}) differs only through the additional term 
$\|GBz\|_1$, which penalizes the inequality constraints and promote the minimality condition~(\ref{condition:minimal}).
}


\section{Experiments}
\label{sec:experiments}

\subsection{Setup}
\subsubsection{Datasets} 
\black{
We evaluated the proposed method on two established benchmark libraries: (i) QPLIB~\cite{qplib}, a comprehensive collection of quadratic programming instances, and (ii) MINLPLib~\cite{bussieck2003minlplib}, a standard repository of mixed-integer nonlinear programming models. 
Maximization problems were reformulated as minimization problems by negating the objective function.
To construct a representative and meaningful testbed, we applied the following selection criteria: (i) the instance conforms to the structure of Model~(\ref{eq:INLP}) or (\ref{eq:NLIP_ineq_slack}), namely linear constraints with integer variables; (ii) the instance is nontrivial, excluding cases solvable to optimality within 10 seconds by the commercial solvers used in our experiments; and (iii) the instance is of moderate size, restricted to at most 700 variables for QPLIB and 500 variables for MINLPLib. This procedure resulted in 53 instances.
The resulting testbed captures a diverse and challenging class of nonlinear integer programs, spanning problem sizes from tens to several hundred variables, varying constraint densities, and both convex and nonconvex objective functions.
}

\subsubsection{Baselines} 
In the experiments, we denote our algorithm by MAPE and compare it with state-of-the-art nonlinear integer programming solvers, including 
(i)~Gurobi 12.0.0~\cite{Gurobi}; 
(ii)~CPLEX 22.1.0~\cite{cplex}; and (iii)~SCIP 9.1.1~\cite{scip}.
\black{We also consider advanced heuristics for integer quadratic programs (iv)~LS-IQCQP~\cite{he2024local} and its successor (v)~QSeek.}
\black{Each baseline solver was run with a time limit of \num{3600} seconds and granted access to all available threads of the testing machine.}

For MAPE, we perform  augmentations from \black{\num{200}} starting points and conducts \black{\num{200000}} parallel extractions.
The hyper-parameters $\lambda_1,\lambda_2$ and $\lambda_3$ are set to \num{0.85}, \num{1} and \num{0.1}, respectively. 
The optimization of unconstrained problems are conducted via Adam with a learning-rate of \num{0.003} (with a PyTorch implementation on GPUs).
Our code is available at \url{https://github.com/NetSysOpt/MAPE}.

\subsubsection{Evaluation Configuration}
All our experiments were conducted on an Nvidia GeForce RTX 3090 GPU and an 12th Gen Intel(R) Core(TM) i9-12900K CPU \black{with 24 threads}, using Python 3.11.5 and PyTorch 2.7.0.

\subsection{Results and Analysis}
In this section, we compare our proposed MAPE with advanced solvers. The results are shown in Table~\ref{table:results}.

\begin{table}
\renewcommand\arraystretch{1.8}
\caption{Computational results on instances from MINLPLib and QPLIB}
\label{table:results}
\centering
\tabcolsep=0.2cm
\resizebox{\textwidth}{!}{
\black{
\begin{tabular}{c@{\hspace{2em}}ccccccccccccccccc}
\toprule
\multirow{2}{*}{}&\multirow{2}{*}{Inst.}& \multirow{2}{*}{\#Vars}& \multirow{2}{*}{\#Cons.} & \multirow{2}{*}{BKS.} & \multicolumn{3}{c}{MAPE} & \multicolumn{2}{c}{Gurobi}& \multicolumn{2}{c}{CPLEX} & \multicolumn{2}{c}{SCIP}& \multicolumn{2}{c}{QSeek} & \multicolumn{2}{c}{LS-IQCQP}  \\ 
\cmidrule(lr){6-8}\cmidrule(lr){9-10}\cmidrule(lr){11-12}\cmidrule(lr){13-14}\cmidrule(lr){15-16}\cmidrule(lr){17-18}
& &  & & & MA & PE & Obj & Time~(Obj.) & Win & Time~(Obj.)&Win & Time~(Obj.) &Win & Time~(Obj.) &Win  & Time~(Obj.) &Win  \\
\cmidrule(lr){1-5}\cmidrule(lr){6-8}\cmidrule(lr){9-10}\cmidrule(lr){11-12}\cmidrule(lr){13-14}\cmidrule(lr){15-16}\cmidrule(lr){17-18}
\ldelim\{{17}{6em}[\text{MINLPLib}]& cardqp\_inlp & 50 & 1 & 3760.7151 & 6.2 & 3.0 & 3760.7151 & 0.0 & $\times$ & 2622.2 &  $\checkmark$ & 795.9 &  $\checkmark$ & 3600.0 &  -- & (4939.6706) & $\checkmark$ \\
&cardqp\_iqp & 50 & 1 & 3760.7151 & 6.0 & 2.9 & 3760.7151 & 0.0 & $\times$ & 2614.2 &  $\checkmark$ & 799.6 &  $\checkmark$ & 3600.0 &  -- & (4939.6706) & $\checkmark$ \\
&graphpart\_clique-50 & 150 & 50 & 2312.0 & 8.4 & 6.7 & 2312.0 & 0.1 & $\times$ & 21.4 &  $\checkmark$ & 161.4 &  $\checkmark$ & 3600.0 &  -- & (4923.0) & $\checkmark$ \\
&graphpart\_clique-60 & 180 & 60 & 3990.0 & 15.0 & 8.2 & 3990.0 & 0.0 & $\times$ & 471.1 &  $\checkmark$ & 1386.9 &  $\checkmark$ & 3600.0 &  -- & (8583.0) & $\checkmark$ \\
&graphpart\_clique-70 & 210 & 70 & 6348.0 & 21.2 & 9.9 & 6348.0 & 0.0 & $\times$ & (6633.0) &  $\checkmark$ & (6854.0) &  $\checkmark$ & 3600.0 &  -- & (7763.0) & $\checkmark$ \\
&qap & 225 & 30 & 388214.0 & 11.3 & 19.9 & 388250.0 & (393734.0) & $\checkmark$ & (410792.0) &  $\checkmark$ & (400948.0) &  $\checkmark$ & 3600.0 &  $\times$ & (475476.0) & $\checkmark$ \\
\ldelim\{{47}{1em}[\text{QPLIB}]&color\_lab2\_4x0 (3913) & 300 & 61 & 42.925 & 42.0 & 15.5 & 43.7 & 4.4 & $\times$ & (50.364) &  $\checkmark$ & (52.125) &  $\checkmark$ & (73.9989) &  $\checkmark$ & (75.9951) & $\checkmark$ \\
&color\_lab3\_3x0 (3931) & 316 & 80 & 79.93 & 42.6 & 15.7 & 80.87 & 5.7 & $\times$ & (80.098) &  $\times$ & (81.166) &  $\checkmark$ & (89.0845) &  $\checkmark$ & (85.49) & $\checkmark$ \\
&color\_lab3\_4x0 (3923)& 395 & 80 & 64.1 & 62.8 & 20.8 & 65.4 & 97.4 & $\times$ & (65.475) &  $\checkmark$ & (67.0) &  $\checkmark$ & (76.8) &  $\checkmark$ & (74.675) & $\checkmark$ \\
&color\_lab6b\_4x20 (3980) & 235 & 48 & 6.325 & 33.6 & 12.0 & 6.325 & 0.4 & $\times$ & (18.45) &  $\checkmark$ & (13.375) &  $\checkmark$ & (8.8) &  $\checkmark$ & (11.7) & $\checkmark$ \\
&qspp\_0\_10\_0\_1\_10\_1 (7139) & 180 & 100 & 621.0 & 0.8 & 6.2 & 621.0 & 0.7 & $\times$ & 37.9 &  $\checkmark$ & 14.0 &  $\checkmark$ & 3600.0 &  -- & 3600.0 & -- \\
&qspp\_0\_11\_0\_1\_10\_1 (7144) & 220 & 121 & 813.0 & 1.3 & 7.8 & 813.0 & 0.4 & $\times$ & 133.4 &  $\checkmark$ & 140.0 &  $\checkmark$ & 3600.0 &  -- & (844.0) & $\checkmark$ \\
&qspp\_0\_12\_0\_1\_10\_1 (7149) & 264 & 144 & 959.0 & 4.2 & 9.4 & 959.0 & 2.4 & $\times$ & 345.7 &  $\checkmark$ & 124.3 &  $\checkmark$ & 3600.0 &  -- & (1002.0) & $\checkmark$ \\
&qspp\_0\_13\_0\_1\_10\_1 (7154) & 312 & 169 & 1159.0 & 5.9 & 11.5 & 1159.0 & 5.5 & $\times$ & 2298.4 &  $\checkmark$ & 1748.0 &  $\checkmark$ & 3600.0 &  -- & (1193.0) & $\checkmark$ \\
&qspp\_0\_14\_0\_1\_10\_1 (7159) & 364 & 196 & 1363.0 & 4.1 & 13.5 & 1363.0 & 1.1 & $\times$ & (1376.0) &  $\checkmark$ & 3301.9 &  $\checkmark$ & (1368.0) &  $\checkmark$ & (1420.0) & $\checkmark$ \\
&qspp\_0\_15\_0\_1\_10\_1 (7164) & 420 & 225 & 1551.0 & 1.5 & 16.0 & 1551.0 & 581.6 & $\checkmark$ & (1625.0) &  $\checkmark$ & (1587.0) &  $\checkmark$ & (1588.0) &  $\checkmark$ & (1657.0) & $\checkmark$ \\
&celar6-sub0 (6324) & 640 & 16 & 159.0 & 88.7 & 66.5 & 159.0 & 5.5 & $\times$ & 297.6 &  $\checkmark$ & 467.8 &  $\checkmark$ & 3600.0 &  -- & (12479.0) & $\checkmark$ \\
&0633 & 75 & 1 & 79.56 & 12.9 & 7.5 & 79.56 & 0.1 & $\times$ & (80.94) &  $\checkmark$ & (81.9572) &  $\checkmark$ & 3600.0 &  $\checkmark$ & (105.7377) & $\checkmark$ \\
&2492 & 196 & 28 & 2724.0 & 13.7 & 16.9 & 2724.0 & (2744.0) & $\checkmark$ & (2766.0) &  $\checkmark$ & (2754.0) &  $\checkmark$ & (2728.0) &  $\checkmark$ & (3148.0) & $\checkmark$ \\
&2512 & 100 & 20 & 135028.0 & 4.0 & 8.2 & 135028.0 & 0.3 & $\times$ & 64.9 &  $\checkmark$ & 170.2 &  $\checkmark$ & 3600.0 &  -- & (175040.0) & $\checkmark$ \\
&2733 & 324 & 36 & 5358.0 & 29.8 & 29.2 & 5358.0 & 59.8 & $\checkmark$ & (5492.0) &  $\checkmark$ & (5414.0) &  $\checkmark$ & (5410.0) &  $\checkmark$ & (6024.0) & $\checkmark$ \\
&2880 & 625 & 50 & 1172056.0 & 62.5 & 60.8 & 1245100.0 & (1202124.0) & $\times$ & (1280184.0) &  $\checkmark$ & (1282352.0) &  $\checkmark$ & (1190140.0) &  $\times$ & (1457718.0) & $\checkmark$ \\
&2957 & 484 & 44 & 3596.0 & 61.8 & 45.3 & 3692.0 & 197.4 & $\times$ & (3962.0) &  $\checkmark$ & (3736.0) &  $\checkmark$ & (3720.0) &  $\checkmark$ & (4898.0) & $\checkmark$ \\
&3307 & 256 & 32 & 1240.0 & 17.7 & 22.3 & 1240.0 & 450.4 & $\checkmark$ & (1336.0) &  $\checkmark$ & (1330.0) &  $\checkmark$ & (1244.0) &  $\checkmark$ & (1658.0) & $\checkmark$ \\
&3347 & 676 & 52 & 3818879.0 & 58.7 & 66.5 & 3834454.0 & (3825877.0) & $\times$ & (3917242.0) &  $\checkmark$ & (3922902.0) &  $\checkmark$ & (3889066.0) &  $\checkmark$ & (4151619.0) & $\checkmark$ \\
&3402 & 144 & 24 & 224416.0 & 6.2 & 12.1 & 224416.0 & 59.4 & $\checkmark$ & (234312.0) &  $\checkmark$ & 2045.8 &  $\checkmark$ & 3600.0 &  -- & (313342.0) & $\checkmark$ \\
&3413 & 400 & 40 & 2192.0 & 39.2 & 36.6 & 2628.0 & 12.0 & $\times$ & (2232.0) &  $\times$ & (2344.0) &  $\times$ & (2232.0) &  $\times$ & (8118.0) & $\checkmark$ \\
&3703 & 225 & 30 & 388214.0 & 15.3 & 19.9 & 388870.0 & (391170.0) & $\checkmark$ & (399684.0) &  $\checkmark$ & (401622.0) &  $\checkmark$ & 3600.0 &  $\times$ & (464698.0) & $\checkmark$ \\
&3750 & 210 & 70 & 6348.0 & 21.7 & 16.9 & 6348.0 & 0.0 & $\times$ & (6477.0) &  $\checkmark$ & (6637.0) &  $\checkmark$ & 3600.0 &  -- & (7763.0) & $\checkmark$ \\
&3751 & 150 & 50 & 2312.0 & 9.1 & 11.3 & 2312.0 & 0.0 & $\times$ & 141.9 &  $\checkmark$ & 359.1 &  $\checkmark$ & 3600.0 &  -- & (4923.0) & $\checkmark$ \\
&3775 & 180 & 60 & 3990.0 & 14.5 & 13.7 & 3990.0 & 0.0 & $\times$ & 718.4 &  $\checkmark$ & 993.0 &  $\checkmark$ & 3600.0 &  -- & (8583.0) & $\checkmark$ \\
&3834 & 50 & 1 & 3760.7151 & 6.6 & 5.1 & 3760.7151 & 0.0 & $\times$ & 2789.2 &  $\checkmark$ & 788.9 &  $\checkmark$ & 3600.0 &  -- & (4939.6706) & $\checkmark$ \\
&6487 & 618 & 309 & 344592.0 & 1.9 & 41.5 & 344592.0 & (345008.0) & $\checkmark$ & (361290.0) &  $\checkmark$ & (348725.0) &  $\checkmark$ & 3600.0 &  -- & (391681.0) & $\checkmark$ \\
&6597 & 600 & 60 & 6491721.0 & 273.3 & 63.5 & 6491721.0 & (6826223.0) & $\checkmark$ & (8022129.0) &  $\checkmark$ & (6749458.0) &  $\checkmark$ & 3600.0 &  -- & (8670743.0) & $\checkmark$ \\
&6647 & 627 & 33 & 2.0 & 101.5 & 65.2 & 3.0 & 1.7 & $\times$ & 18.2 &  $\times$ & 128.7 &  $\times$ & 3600.0 &  $\times$ & (22.0) & $\checkmark$ \\
&0752 & 250 & 1 & -24071.0 & 27.8 & 15.7 & -24071.0 & 12.4 & $\times$ & (-23987.0) & $\checkmark$ & (-23757.0) &  $\checkmark$ & 3600.0 &  -- & 3600.0 &-- \\
&2315 & 595 & 13090 & -29432.0 & 4.5 & 195.4 & -9740.0 & (-26935.0) & $\times$ & (-18454.0) & $\times$ & (-14940.0) &  $\times$ & 3600.0 &  $\times$ & (-26935.0) & $\times$ \\
&2357 & 240 & 2240 & -647.0 & 24.9 & 28.1 & -647.0 & 7.0 & $\times$ & 244.7 & $\checkmark$ & 703.5 &  $\checkmark$ & 3600.0 &  -- & 3600.0 &-- \\
&2359 & 306 & 3264 & -648.0 & 32.4 & 43.3 & -648.0 & 5.3 & $\times$ & 87.4 & $\checkmark$ & 486.9 &  $\checkmark$ & 3600.0 &  -- & (-638.0) & $\checkmark$ \\
&3584 & 528 & 10912 & -25386.0 & 18.9 & 157.6 & -23435.0 & (-24020.0) & $\times$ & (-18038.0) & $\checkmark$ & (-12510.0) &  $\checkmark$ & 3600.0 &  $\times$ & (-22989.0) & $\checkmark$ \\
&3752 & 462 & 6160 & -1306.0 & 55.4 & 89.0 & -1297.0 & 864.5 & $\times$ & (-1266.0) & $\checkmark$ & (-1274.0) &  $\checkmark$ & (-1250.0) &  $\checkmark$ & (-1256.0) & $\checkmark$ \\
&3757 & 552 & 8096 & -563.0 & 59.0 & 128.9 & -521.0 & 21.4 & $\times$ & 77.7 & $\times$ & 552.5 &  $\times$ & 3600.0 &  $\times$ & 3600.0 & $\times$ \\
&3762 & 90 & 480 & -296.0 & 4.9 & 6.9 & -296.0 & 2.4 & $\times$ & 25.6 & $\checkmark$ & 62.2 &  $\checkmark$ & 3600.0 &  -- & (-272.0) & $\checkmark$ \\
&3772 & 380 & 4560 & -940.0 & 42.3 & 60.0 & -940.0 & 94.1 & $\times$ & 2599.3 & $\checkmark$ & 3151.8 &  $\checkmark$ & (-834.0) &  $\checkmark$ & (-936.0) & $\checkmark$ \\
&3803 & 190 & 2280 & -7360.0 & 18.8 & 23.5 & -7360.0 & 4.5 & $\times$ & 544.0 & $\checkmark$ & 887.9 &  $\checkmark$ & 3600.0 &  -- & (-6980.0) & $\checkmark$ \\
&3841 & 300 & 4600 & -1817.0 & 24.1 & 51.3 & -1799.0 & 1235.9 & $\times$ & (-1627.0) & $\checkmark$ & (-1623.0) &  $\checkmark$ & 3600.0 &  $\times$ & (-1769.0) & $\checkmark$ \\
&3860 & 435 & 8120 & -20161.0 & 10.4 & 105.8 & -19756.0 & (-18987.0) & $\checkmark$ & (-14767.0) & $\checkmark$ & (-16191.845) &  $\checkmark$ & 3600.0 &  $\times$ & 3600.0 & $\times$ \\
&3883 & 182 & 1456 & -788.0 & 16.8 & 18.5 & -788.0 & 338.8 & $\checkmark$ & 3251.2 & $\checkmark$ & (-776.0) &  $\checkmark$ & (-787.0) &  $\checkmark$ & (-777.0) & $\checkmark$ \\
&5935 & 100 & 1237 & -4758.0 & 5.6 & 6.4 & -4758.0 & 7.9 & $\times$ & 468.9 & $\checkmark$ & 907.4 &  $\checkmark$ & (-4164.0) &  $\checkmark$ & 3600.0 & -- \\
&5944 & 100 & 2475 & -1829.0 & 8.5 & 19.8 & -1829.0 & 0.8 & $\times$ & 58.7 & $\checkmark$ & 67.9 &  $\checkmark$ & (-1662.0) &  $\checkmark$ & 3600.0 & -- \\
&5962 & 150 & 2793 & -6962.0 & 39.4 & 59.7 & -6962.0 & 52.7 & $\times$ & (-6318.0) & $\checkmark$ & (-6296.0) &  $\checkmark$ & (-4369.0) &  $\checkmark$ & 3600.0 & -- \\
&5971 & 150 & 5587 & -2377.0 & 65.9 & 46.9 & -2377.0 & 2.9 & $\times$ & 275.6 & $\checkmark$ & 149.3 &  $\checkmark$ & (-1719.0) &  $\checkmark$ & 3600.0 & -- \\
&5980 & 150 & 8381 & -895.0 & 14.8 & 58.5 & -895.0 & 38.1 & $\times$ & 65.9 & $\times$ & 69.8 &  $\times$ & (-875.0) &  $\checkmark$ & 3600.0 & -- \\
\midrule
& & && & & & & & 11/53& &47/53 & &48/53 & &20/30 & & 42/45\\

\bottomrule
\end{tabular}
}
}
\end{table}

\subsubsection{Metric}
The columns labeled \enquote{Inst.}, \enquote{\# Vars.}, \enquote{\# Cons.} and \enquote{Bkv.} display, respectively, the instance name, the number of variables and constraints, and the best-known objective value for each instance. For MAPE, we report three key metrics: (i) computational time for the multi-start augmentation (\enquote{MA}) procedure, which involves computing initial solutions, (ii) computational time for the parallel extraction (\enquote{PE}), and (iii) the best objective achieved (\enquote{Obj}).
The last ten columns present the performance (\enquote{Time~(Obj.)}) for each solver, alongside a comparison with MAPE. For \enquote{Time~(Obj.)}, values without curly brackets represent the computational time taken by the solver to find a solution with an objective value equal to \enquote{Bks.} (but not necessarily with proven optimality), while \enquote{$(\cdot)$} indicates timeouts, with bracketed values denoting the best objectives achieved in those cases. Moreover, \enquote{Win} indicates whether MAPE outperforms the baseline solver in terms of both solution quality and computational time. 
\black{This comparison is not applicable (denoted as \enquote{--}) when both MAPE and QSeek (LS-IQCQP) reach the best known solution, since detailed runtime data for the latter are unavailable.}

\subsubsection{Analysis}
The computational results are summarized in Table~\ref{table:results}.
MAPE outperforms CPLEX and SCIP on 90\% of the instances, and \black{surpasses QSeek and LS-IQCQP on more than 66\% and 93\% of comparable cases, respectively}.
While Gurobi attains optimal or best-known solutions for more than half of the instances within seconds, MAPE still delivers superior performance on 20\% of the instances and prominently surpasses all baseline solvers on instances \enquote{2492}, \enquote{2733}, \enquote{3307}, \enquote{3883} and \enquote{7164}. 
\black{Notably, whereas these general-purpose solvers incorporate sophisticated components such as presolve routines and primal heuristics, MAPE relies on a lightweight implementation of only a few hundred lines of Python code while still producing high-quality solutions. Its GPU-based design also suggests potential complementarity with CPU-based solvers such as Gurobi. Overall, MAPE provides a competitive alternative to state-of-the-art solvers.}


\section{Conclusion}
\label{sec:conclusion}
In this paper, we propose MAPE, a heuristic approach designed to extract the Graver basis elements and solve nonlinear integer programs through a multi-start augmentation procedure. 
Unlike traditional approaches that struggle with sequential processing, MAPE excels in parallelism, enabling faster execution. 
Experimental results on public datasets demonstrate MAPE's ability to swiftly deliver high-quality solutions, outperforming state-of-the-art solvers. Finally, MAPE does not rely on any general-purpose solver or sophisticated presolving procedures, and \black{could be massively paralleled on GPUs}. 

\section*{Acknowledgements}
This work was supported by the National Key R\&D Program of China (Grant No.~\seqsplit{2022YFA1003900}). Akang Wang also acknowledges support from National Natural Science Foundation of China (Grant No.~\seqsplit{12301416}), Guangdong Basic and Applied Basic Research Foundation (Grant No.~\seqsplit{2024A1515010306}), Shenzhen Science and Technology Program (Grant No.~\seqsplit{JCYJ20250604191330040}) and Hetao Shenzhen-Hong Kong Science and Technology Innovation Cooperation Zone Project (No. HZQSWS-KCCYB-2024016).

\bibliographystyle{elsarticle-num} 
\bibliography{ref}

\newpage
\appendix
\section{Additional Preliminaries}
\label{apdx:lll}
Given a lattice $\mathcal L(B)$, then $B=(B_1,...,B_d)$ is a basis but not unique, and we aim to identify an appropriate basis that is nearly orthogonal with short base vectors.
While the vectors $\{B_i^*\}_{i=1}^d$ obtained from Gram-Schmidt orthogonalization (GSO) $B_i^*=B_i-\sum\limits_{j<i}\frac{\langle B_i,B_j^*\rangle}{\langle B_j^*,B_J^*\rangle}B_j^*$ enjoy favorable properties, they do not form a basis.
Alternatively, we define the \textit{reduced basis}.

\textbf{Definition:}
A basis \( B = (B_1, \dots, B_d) \) of a lattice \( \mathcal{L} \) is called \textit{reduced} if it satisfies the following conditions, where \( B^* \) denotes the corresponding Gram-Schmidt vectors:
\begin{itemize}
    \item[(i)] \( \left| \frac{\langle B_i, B_j^* \rangle}{\langle B_j^*, B_j^* \rangle} \right| \leq \frac{1}{2} \);
    \item[(ii)] \( \Vert B_j^* \Vert^2 \leq 2 \Vert B_{j+1}^* \Vert^2 \).
\end{itemize}

Algorithm~\ref{alg:LLL} presents the classic basis reduction algorithm proposed by~\cite{lenstra1982factoring}, which is known as the \textit{LLL algorithm} and runs in polynomial time~\cite{schrijver1998theory}.
For the sake of simplicity, let $\lceil \mu\rfloor$ denote the integer nearest to $\mu$.
\begin{algorithm}
\caption{The LLL algorithm}
\label{alg:LLL}
\begin{algorithmic}[1] 
\Statex \textbf{INPUT:} a basis $B=(B_1,...,B_d)$ of the lattice $\mathcal L$
\Statex \textbf{OUTPUT:} the reduced basis $B$ of $\mathcal L$
\State $B_i^*=B_i-\sum\limits_{j<i}\mu_{ij}B_j^*$ with $\mu_{ij}=\frac{\langle B_i,B_j^*\rangle}{\langle B_j^*,B_J^*\rangle}$ \hspace{1 in}\#compute GSO
\State $i\gets 2$
\While{$i\leq d$}
\For {$j=i-1, i-2, ...,1$}
\State $B_i\gets B_i- \lceil\mu_{ij}\rfloor B_j$, update GSO
\If {$i>1$ and $\Vert B_{i-1}^*\Vert ^2>2\Vert B_i^*\Vert ^2$}
\State Exchange $B_i,B_{i-1}$, update GSO, $i\gets i-1$
\Else
\State $i\gets i+1$
\State \textbf{return} $B$
\EndIf
\EndFor
\EndWhile
\end{algorithmic}
\end{algorithm}

\end{document}